\RequirePackage{fix-cm}
\documentclass[smallextended]{svjour3}       
\smartqed  
\usepackage{graphicx}
\usepackage{epsfig,epsfig}
\usepackage{amsmath}
\usepackage{amssymb}
\usepackage{booktabs,arydshln,multirow}
\usepackage{setspace}
\newtheorem{thm}{Theorem}
\newtheorem{lem}{Lemma}

\newtheorem{defn}{Definition}
\newtheorem{rem}{Remark}

\newcommand{\be}{\begin{eqnarray}}
\newcommand{\ee}{\end{eqnarray}}
\begin{document}

\title{A two-grid finite element approximation for a nonlinear time-fractional Cable equation
}

\titlerunning{Two-grid FE approximation for nonlinear time-fractional Cable equation}        

\author{Yang Liu   \and Yan-Wei Du \and  Hong Li  \and Jin-Feng Wang
}


\institute{Y. Liu \and Y.W. Du\and H. Li\and J.F. Wang\at
              School of Mathematical Sciences, Inner Mongolia University, Hohhot, 010021, China \\
              Tel.: +86-471-4994853\\
              Fax: +86-471-4991650\\
              \email{mathliuyang@aliyun.com}           
           \and
           J.F. Wang  \at
              School of Statistics and Mathematics, Inner Mongolia University of Finance and Economics, Hohhot, 010070, China
}

\date{Received: date / Accepted: date}

\maketitle

\begin{abstract}
In this article, a nonlinear fractional Cable equation is solved by a two-grid algorithm combined with finite element (FE) method. A temporal second-order fully discrete two-grid FE scheme, in which the spatial direction is approximated by two-grid FE method and the integer and fractional derivatives in time are discretized by second-order two-step backward difference method and second-order weighted and shifted Gr\"{u}nwald difference (WSGD) scheme, is presented to solve nonlinear fractional Cable equation.
The studied algorithm in this paper mainly covers two steps: First, the numerical solution of nonlinear FE scheme on the coarse grid is solved; Second, based on the solution of initial iteration on the coarse grid, the linearized FE system on the fine grid is solved by using Newton iteration. Here, the stability based on fully discrete two-grid method is derived. Moreover, the a priori estimates with second-order convergence rate in time is proved in detail, which is higher than the L1-approximation result with $O(\tau^{2-\alpha}+\tau^{2-\beta})$. Finally, the numerical results by using the two-grid method and FE method are calculated, respectively, and the CPU-time is compared to verify our theoretical results.
\end{abstract}

\section{Introduction}
\label{sec:1}
Fractional partial differential equations (FPDEs), which have a lot of applications in the realm of science, mainly include space FPDEs \cite{Wanghyd,Me2,Sousa2,Huangcm1,Xiaoa1,Xuqw1,Roop},
time FPDEs \cite{Baleanu5,Caojx2,Dinghf1,Xucj1,Heyn1,Sunhw1,Atangana2,Mclean1,Chencm4,Lukashchuk,Sunzzwuxn} and space-time FPDEs \cite{Qianqian1,Liufw2,El-Wakil}. The construction of numerical methods for FPDEs has attracted great attention of many scholars. For example, finite element (FE) methods have been successfully applied to solving many FPDEs in the current literatures. In \cite{Roop}, Roop gave FE method for fractional advection dispersion equations on bounded domains in two-dimensions. In \cite{Liufw2}, Feng \emph{et al.} studied FE method for diffusion equation with space-time fractional derivatives. In \cite{MA1}, Ma \emph{et al.} used moving FE methods to solve
space fractional differential equations. Li \emph{et al.} in \cite{Ljc2} gave some numerical theories on FE methods for Maxwell's
equations. In \cite{Liuyli} Liu et al. proposed a mixed FE method for a fourth-order time-FPDE with first-order convergence rate in time.
In \cite{Liuylh2}, Liu \emph{et al.} solved a time-fractional reaction-diffusion problem with fourth-order space derivative term by using FE method and L1-approximation. In \cite{Jinbt6} Jin \emph{et al.} used a FE method to solve the space-fractional parabolic equation, and gave the error analysis.
In \cite{Zengfh1} Zeng \emph{et al.} used FE approaches combined with finite difference method for solving the time fractional
subdiffusion equation. In \cite{Ford} Ford \emph{et al.} studied a FE method for time FPDEs and obtained optimal convergence error estimates.
Bu \emph{et al.} \cite{Tang1} discussed Galerkin FE method for Riesz space fractional diffusion equations in two-dimensional case.
In \cite{LCP1}, Li \emph{et al.} applied FE method to solving nonlinear fractional subdiffusion and superdiffusion equations.
In \cite{Dengwh7}, Deng solved fractional Fokker-Planck equation with space and time derivatives by using FE method.
In \cite{Wu1}, Zhang \emph{et al.} implemented FE method for solving a modified fractional diffusion equation in two-dimensional case.
\par
In this article, we will consider a two-grid FE algorithm for solving a nonlinear time-fractional Cable equation
\begin{equation}\label{0.1}
\frac{\partial u}{\partial t}=-^R_0\partial_{t}^\alpha u+^R_0\partial_{t}^\beta\Delta u-\mathcal{F}(u)+g(\textbf{x},t),(\textbf{x},t)\in\Omega\times J,
\end{equation}
which covers boundary condition
 \begin{equation}\label{0.2}
u(\textbf{x},t)= 0,  (\textbf{x},t)\in\partial\Omega\times \bar{J},
\end{equation}
and initial condition
\begin{equation}\label{0.3}
u(\textbf{x},0)=u_0(\textbf{x}), \textbf{x}\in\Omega,
\end{equation}
where $\Omega$ is a bounded convex polygonal sub-domain of $R^d(d\leq 2)$, whose  boundary $\partial\Omega$ is $Lipschitz$ continuous. $J=(0,T]$ is the time interval with the upper bound $T$. The source item $g(\textbf{x},t)$ and the initial function $u_0(\textbf{x})$ are given known functions. For the nonlinear term $\mathcal{F}(u)$, there exists a constant $C>0$ such that $|\mathcal{F}(u)| \leq C|u|$ and $|\mathcal{F}'(u)| \leq C$. And $^R_0\partial_{t}^\gamma w(\textbf{x},t)$ is Riemann-Liouville fractional derivative with $\gamma\in (0,1)$ given in Definition \ref{y12}.
\par
The fractional Cable equation \cite{Henry,Bisquert,Henry1,Tomovski1}, which reflects the anomalous electro-diffusion in nerve cells, is an important mathematical model. For the fractional Cable equation, we can find some numerical methods, such as finite difference methods \cite{Chencm1,Zhanglm1,Yubo,Yusteq,QQ}, orthogonal spline collocation method \cite{Zhanghx1}, spectral approximations \cite{Bhrawy1,Xucj6}, finite element methods \cite{Zhuangph1,Liuylh}. Chen \emph{et al}. \cite{Chencm1}, Hu and Zhang \cite{Zhanglm1}, Yu and Jiang \cite{Yubo}, Quintana-Murillo and Yuste \cite{Yusteq}, presented and analyzed some finite difference schemes to numerically solve the fractional cable equation from different perspectives. Liu \emph{et al}. \cite{QQ}, solved numerically the fractional cable equation by using two implicit numerical schemes. Zhang \emph{et al}. \cite{Zhanghx1} proposed and analyzed the discrete-time orthogonal spline collocation method for the two-dimensional case of fractional cable equation.
Bhrawy and Zaky \cite{Bhrawy1} presented a Jacobi spectral collocation approximation for numerically solving nonlinear two-dimensional
fractional cable equation covering Caputo fractional
derivative. Lin \emph{et al}. \cite{Xucj6} developed spectral approximations combined with finite difference method for looking for the numerical solution of
the fractional Cable equation. Recently, Zhuang \emph{et al}. \cite{Zhuangph1}, Liu \emph{et al}. \cite{Liuylh} studied and analyzed Galerkin finite element methods for the fractional Cable equation with Riemann-Liouville derivative, respectively, and did some different analysis based on different approximate formula for fractional derivative. Here, we will consider a two-grid FE algorithm combined with a higher-order time approximation to seek the numerical solutions of nonlinear fractional Cable equation.
\par
Two-grid FE algorithm was presented and developed by Xu \cite{Xuj1,Xuj2}. Owing to holding the advantage of saving computation time, many computational scholars have used well the method to numerically solve integer-order partial differential equations(such as Dawson and Wheeler\cite{Dawson} for nonlinear parabolic equations; Zhong et al. \cite{zhonglq} for for time-harmonic Maxwell equations; Mu and Xu \cite{Xuj4} for mixed Stokes-Darcy model; Chen and Chen \cite{Cheny3} for nonlinear reaction-diffusion equations; Bajpai and Nataraj \cite{Baj} for Kelvin-Voigt model; Wang \cite{Wangws1} for semilinear evolution equations with
positive memory) and developed some new numerical techniques based on the idea of two-grid algorithm (two-grid expanded Mixed FE methods in Chen \emph{et al}. \cite{Cheny1}, Wu and Allen \cite{Wul1}, Liu \emph{et al}. \cite{Ruih3}; two-grid finite volume element method in Chen and Liu \cite{Ruih2}). Until recently, in \cite{Liuylh1}, the two-grid FE method was presented to solve the nonlinear fourth-order fractional differential equations with Caputo fractional derivative. However, the Caputo time fractional derivative was approximated by L1-formula and the only ($2-\alpha$)-order convergence rate in time was arrived at in \cite{Liuylh1}.
\par
In this article, our main task is to look for the numerical solution of nonlinear fractional Cable equation (\ref{0.1}) with initial and boundary condition by using two-grid FE method with higher-order time approximate scheme \cite{Dengwhtw,Vongwang1} and to discuss the numerical theories on stability and a priori estimate analysis for this method. In \cite{Dengwhtw}, Tian \emph{et al.} approximated the Riemann-Liouville fractional derivative by proposing a new higher-order WSGD operator, then discussed some finite difference scheme based on this operator. Considering this idea of WSGD operator, Wang and Vong \cite{Vongwang1} presented the compact difference scheme for the modified anomalous sub-diffusion equation with $\alpha$-order Caputo fractional derivative, in which the Caputo fractional derivative covering order $\alpha\in (0,1)$ is approximated by applying the idea of WSGD operator, and an extension
for this idea was also made to discuss a compact difference scheme for the fractional diffusion-wave equation. However, the theories of the FE methods based on
the idea of WSGD operator have not been studied and discussed. Especially, the two-grid FE algorithm combined with the idea of the WSGD operator has not been reported in the current literatures. Here, we will study the two-grid FE scheme with WSGD operator for solving nonlinear fractional Cable equation, derive the stability of the studied method, and prove a priori estimate results with second-order convergence rate, which is higher than time convergence rate $O(\tau^{2-\alpha}+\tau^{2-\beta})$ obtained by usual L1-approximation. Finally, we do some numerical computations by using the current method and standard nonlinear FE method, respectively, and find that our method in CPU-time is more efficient than standard nonlinear FE method.
\par
Throughout this article, we will denote $C>0$ as a constant, which is
free of the spatial coarse grid size $H$, fine step length $h$, and time mesh size $\tau$. Further, we define the natural inner product in $L^2(\Omega)$ or $(L^2(\Omega))^2$ by
~$(\cdot,\cdot)$ equipped with norm $\|\cdot\|$.
\par
The remaining outline of the article is as follows. In section 2, some definitions of fractional derivatives, lemmas on time approximations and two-grid algorithm combined with second-order scheme in time are given. In section 3, the analysis of stability of two-grid FE method is made. In section 4, a priori errors of two-grid FE algorithm are proved. In section 5, some numerical results by using two-grid FE method and standard nonlinear FE method are computed and some comparisons of computing time are done. Finally, some remarking conclusions on two-grid FE algorithm proposed in this paper are shown in section 6.
\section{Fractional derivatives and two-grid FE method}
\label{sec:2}
\subsection{Fractional derivatives and approximate formula}
\label{subsec:2.1}
In many literatures, we can get the following definitions on fractional derivatives of Caputo type and Riemann-Liouville type. At the same time, we need to give some useful lemmas in the subsequent theoretical analysis.
\begin{defn}\label{y12}£ºThe $\gamma$-order $(0<\gamma<1)$ fractional derivative of Riemann-Liouville type for the function $w(t)$ is defined as
\begin{equation}\label{1.2}
^R_0\partial_{t}^\alpha w(t) =\frac{1}{\Gamma(1-\gamma)}\frac{d}{dt}\int_0^t\frac{w(\tau)}{(t-\tau)^{\gamma}}d\tau.
\end{equation}
\end{defn}
\begin{defn}\label{y1l}£ºThe $\gamma$-order $(0<\gamma<1)$ fractional derivative of Caputo type for the function $w(t)$ is defined by
\begin{equation}\label{1.2}
_0^C\partial_{t}^\gamma w(t) =\frac{1}{\Gamma(1-\gamma)}\int_0^t\frac{w'(\tau)}{(t-\tau)^\gamma}d\tau,
\end{equation}
where $\Gamma(\cdot)$ is Gamma function.
\end{defn}
\begin{lem}\label{lemma1}\cite{Sunzzjc} The relationship between Caputo fractional derivative and Riemann-Liouville fractional derivative can be given by
\begin{equation}\label{1.8}
^R_{0}\partial_{t}^\gamma w(t) =_{0}^{C}\partial_{t}^\gamma w(t)+\frac{w(0)t^{-\gamma}}{\Gamma(1-\gamma)}.
\end{equation}
\end{lem}
\begin{lem}\label{lemma5}
For $0<\gamma<1$, the following approximate formula \cite{Dengwhtw,Vongwang1} with second-order accuracy at time $t=t_{n+1}$ holds
\begin{equation}\begin{split}\label{3.1}
^R_0\partial_{t}^\gamma w(\mathbf{x},t_{n+1})&
=\sum_{i=0}^{n+1}\frac{p_{\gamma}(i)}{\tau^{\gamma}}w(\mathbf{x},t_{n+1-i})+O(\tau^{2}),
\end{split}\end{equation}
where
\begin{equation}\label{3.2}
p_{\gamma}(i)=\begin{cases}
\frac{\gamma+2}{2}g_{0}^{\gamma},&\text{if} \ i=0,\\
\frac{\gamma+2}{2}g_{i}^{\gamma}+\frac{-\gamma}{2}g_{i-1}^{\gamma},&\text{if} \ i>0,
\end{cases}
\end{equation}
\begin{equation}\begin{split}g_{0}^{\gamma}=1,
g_{i}^{\gamma}=\frac{\Gamma(i-\gamma)}{\Gamma(-\gamma)\Gamma(i+1)},
g_{i}^{\gamma}=\Big{(}1-\frac{\gamma+1}{i}\Big{)}g_{i-1}^{\gamma},~i\geq1.
\end{split}\end{equation}
\end{lem}
\begin{lem}\label{lemma8}For series $\{g_{i}^{\gamma}\}$ defined in lemma \ref{lemma5}, we have
\begin{equation}\begin{split}\label{lemma3}
g_{0}^{\gamma}=1>0,~g_{i}^{\gamma}<0, (i=1,2\cdots),\sum_{i=1}^{\infty}g_{i}^{\gamma}=-1,
\end{split}
\end{equation}
\end{lem}
\begin{lem}\label{lemma4} For series $\{p_{\gamma}(i)\}$ given by (\ref{3.2}), the following inequality holds for any integer $n$
\begin{equation}\begin{split}\label{4.21}
\sum_{i=0}^{n+1}|p_{\gamma}(i)|\leq C.
\end{split}
\end{equation}
\end{lem}
\textbf{Proof}. Noting that the notation (\ref{3.2}), we have
\begin{equation}\begin{split}\label{4.4}
&\sum_{i=0}^{n+1}|p_{\gamma}(i)|=\frac{\gamma+2}{2}g_{0}^{\gamma}
+\sum_{i=1}^{n+1}\Big{|}\frac{\gamma+2}{2}g_{i}^{\gamma}+\frac{-\gamma}{2}g_{i-1}^{\gamma}\Big{|}.
\end{split}
\end{equation}
Applying triangle inequality and lemma \ref{lemma8}, we arrive at
\begin{equation}\begin{split}\label{4.41}
\sum_{i=0}^{n+1}|p_{\gamma}(i)|\leq&
\Big{(}\frac{\gamma+2}{2}g_{0}^{\gamma}
+\sum_{i=1}^{n+1}\Big{|}-\frac{\gamma+2}{2}g_{i}^{\gamma}\Big{|}
+\sum_{i=1}^{n+1}\Big{|}\frac{-\gamma}{2}g_{i-1}^{\gamma}\Big{|}\Big{)}
\\\leq& \Big{(}(\gamma+1)g_{0}^{\gamma}
+\frac{\gamma+2}{2}\sum_{i=1}^{n+1}-g_{i}^{\gamma}
+\frac{\gamma}{2}\sum_{i=1}^{n}-g_{i}^{\gamma}\Big{)}
\\\leq& 2\gamma+2.
\end{split}
\end{equation}
So, we get the conclusion of lemma.
\begin{lem}\label{lm1}\cite{Dengwhtw,Vongwang1} Let $\{p_{\gamma}(i)\}$ be defined as in (\ref{3.2}). Then for any positive integer $L$ and real vector $(w^{0}, w^{1},\cdots,w^{L})\in R^{L+1}$, it holds that
\begin{equation}\label{3.41}
\sum_{n=0}^{L}\Big{(}\sum_{i=0}^{n}p_{\gamma}(i)w^{n-i}\Big{)}w^{n}\geq0.
\end{equation}
 \end{lem}
\begin{rem}
Based on the relationship (\ref{1.8}) between Caputo fractional derivative and Riemann-Liouville fractional derivative, we easily find that the equality $^R_{0}\partial_{t}^\gamma w(t) =_{0}^{C}\partial_{t}^\gamma w(t)$ with $w(0)=0$ holds. Further, it is not hard to know that the second-order discrete formula (\ref{3.1}) in lemma \ref{lemma5} can also approximate the Caputo fractional derivative (\ref{1.2}) with zero initial value.
\end{rem}

\subsection{Two-grid algorithm based on FE scheme}
\label{subsec:2.2}
To give the fully discrete analysis, we should approximate both integer and fractional derivatives. The grid points in the time interval $[0,T]$ are labeled as $t_{i}=i\tau$, $i=0,1,2,\ldots,M$, where $\tau=T/M$ is the time step length. We define $w^{n}=w(t_{n})$ for a smooth function on $[0,T]$ and $\delta^{n}_tw^{n}=\frac{w^n-w^{n-1}}{\tau}$.
\par
Using the approximate formula (\ref{3.1}) and two-step backward Euler approximation, then applying Green's formula, we find $u^{n+1}:[0,T]\mapsto
H^1_0$ to arrive at the weak formulation of (\ref{0.1})-(\ref{0.3}) for any $v\in H^1_0$ as
\par
Case $n=0$:
\begin{equation}\begin{split}\label{3.30}
\Big{(}\delta^{1}_t u^{1},v\Big{)}+\sum_{i=0}^{1}\frac{p_{\alpha}(i)}{\tau^{\alpha}}(u^{1-i},v)
+&\sum_{i=0}^{1}\frac{p_{\beta}(i)}{\tau^{\beta}}(\nabla u^{1-i},\nabla v)+(\mathcal{F}(u^{1}),v)\\=&(g^{1},v)+(\bar{e}^1_1,v)+(\bar{e}^1_2,v)+(\Delta\bar{e}^1_3,v),
\end{split}\end{equation}
\par
Case $n\geq1$:
\begin{equation}\begin{split}\label{3.10}
\Big{(}\frac32\delta^{n+1}_tu^{n+1}-&\frac12\delta^{n}_tu^{n},v\Big{)}+\sum_{i=0}^{n+1}\frac{p_{\alpha}(i)}{\tau^{\alpha}}(u^{n+1-i},v)
+\sum_{i=0}^{n+1}\frac{p_{\beta}(i)}{\tau^{\beta}}(\nabla u^{n+1-i},\nabla v)\\+&(\mathcal{F}(u^{n+1}),v)=(g^{n+1},v)+(\bar{e}^{n+1}_1,v)+(\bar{e}^{n+1}_2,v)+(\Delta\bar{e}^{n+1}_3,v),
\end{split}\end{equation}
where
\begin{equation}\begin{split}\label{3.21}
\bar{e}^{n+1}_1=\left\{\begin{aligned}&\delta^{1}_tu^{1}-u(t_1)=O(\tau),&&n=0,
\\&\frac32\delta^{n+1}_tu^{n+1}-\frac12\delta^{n}_tu^{n}-u_t(t_{n+1})=O(\tau^2),&&n\geq1,\end{aligned}
\right.
\end{split}\end{equation}
\begin{equation}\begin{split}\label{3.22}
\bar{e}^{n+1}_2=O(\tau^2),n\geq0,
\end{split}\end{equation}
\begin{equation}\begin{split}\label{3.23}
\bar{e}^{n+1}_3=O(\tau^2),n\geq0.
\end{split}\end{equation}
For formulating finite element algorithm, we choose finite element space $V_h\subset H_0^1$ as
\begin{equation}\begin{split}\label{2.5} &V_h=\{v_h\in  H_0^1(\Omega)\cap C^0(\overline{\Omega})\mid~v_h|_{e}\in Q_m(e), \forall e\in \mathcal{K}_h\},
 \end{split}
 \end{equation}
 where $\mathcal{K}_h$ is the quasiuniform rectangular partition for the spatial domain $\Omega$.
 \par
Then, we find
$u^{n+1}_h\in V_{h}(n=0,1,\cdots,N_{\tau}-1)$ to formulate a standard nonlinear finite element system for any $v_h\in V_h$ as
\par
Case $n=0$:
\begin{equation}\begin{split}\label{3.31}
\Big{(}\delta^{1}_t u^{1}_h,v_h\Big{)}+\sum_{i=0}^{1}\frac{p_{\alpha}(i)}{\tau^{\alpha}}(u_h^{1-i},v_h)&+\sum_{i=0}^{1}\frac{p_{\beta}(i)}{\tau^{\beta}}(\nabla u_h^{1-i},\nabla v_h)\\&+(\mathcal{F}(u_h^{1}),v_h)=(g^{1},v_h),
\end{split}\end{equation}
\par
Case $n\geq 1$:
\begin{equation}\begin{split}\label{3.3}
\Big{(}\frac32\delta^{n+1}_tu_h^{n+1}&-\frac12\delta^{n}_tu^{n}_h,v_h\Big{)}+\sum_{i=0}^{n+1}\frac{p_{\alpha}(i)}{\tau^{\alpha}}(u_h^{n+1-i},v_h)
\\&+\sum_{i=0}^{n+1}\frac{p_{\beta}(i)}{\tau^{\beta}}(\nabla u_h^{n+1-i},\nabla v_h)+(\mathcal{F}(u_h^{n+1}),v_h)=(g^{n+1},v_h).
\end{split}\end{equation}
For improving the finite element discrete system (\ref{3.31})-(\ref{3.3}), we consider the following two-grid FE system based on the coarse grid $\mathfrak{T}_H$ and the fine grid $\mathfrak{T}_h$.
\\
\textbf{Step I:} First, the following nonlinear system based on the coarse grid $\mathfrak{T}_H$ is solved by finding the solution $u^{n+1}_H:[0,T]\mapsto V_H\subset V_h$ such that
\par
Case $n=0$:
\begin{equation}\begin{split}\label{3.32}
\Big{(}\delta^{1}_t u^{1}_H,v_H\Big{)}+\sum_{i=0}^{1}\frac{p_{\alpha}(i)}{\tau^{\alpha}}(u_H^{1-i},v_H)&+\sum_{i=0}^{1}\frac{p_{\beta}(i)}{\tau^{\beta}}(\nabla u_H^{1-i},\nabla v_H)\\&+(\mathcal{F}(u_H^{1}),v_H)=(g^{1},v_H),
\end{split}\end{equation}
\par
Case $n\geq 1$:
\begin{equation}\begin{split}\label{3.33}
&\Big{(}\frac32\delta^{n+1}_tu_H^{n+1}-\frac12\delta^{n}_tu^{n}_H,v_H\Big{)}+\sum_{i=0}^{n+1}\frac{p_{\alpha}(i)}{\tau^{\alpha}}(u_H^{n+1-i},v_H)
\\+&\sum_{i=0}^{n+1}\frac{p_{\beta}(i)}{\tau^{\beta}}(\nabla u_H^{n+1-i},\nabla v_H)+(\mathcal{F}(u_H^{n+1}),v_H)=(g^{n+1},v_H).
\end{split}\end{equation}
\\
\textbf{Step II:} Second, based on the solution $u^{n+1}_H\in V_H$ on the coarse grid $\mathfrak{T}_H$, the following linear system on the fine grid $\mathfrak{T}_h$, is considered by looking for $U^{n+1}_h:[0,T]\mapsto V_h$ such that
\par
Case $n=0$:
\begin{equation}\begin{split}\label{3.34}
\Big{(}\delta^{1}_t U^{1}_h,v_h\Big{)}&+\sum_{i=0}^{1}\frac{p_{\alpha}(i)}{\tau^{\alpha}}(U_h^{1-i},v_h)+\sum_{i=0}^{1}\frac{p_{\beta}(i)}{\tau^{\beta}}(\nabla U_h^{1-i},\nabla v_h)\\&+(\mathcal{F}(u_{H}^{1})+\mathcal{F}'(u_{H}^{1})(U_{h}^{1}-u_{H}^{1}),v_h)=(g^{1},v_h),
\end{split}\end{equation}
\par
Case $n\geq 1$:
\begin{equation}\begin{split}\label{3.35}
&\Big{(}\frac32\delta^{n+1}_tU_h^{n+1}-\frac12\delta^{n}_tU^{n}_h,v_h\Big{)}\\&+\sum_{i=0}^{n+1}\frac{p_{\alpha}(i)}{\tau^{\alpha}}(U_h^{n+1-i},v_h)
+\sum_{i=0}^{n+1}\frac{p_{\beta}(i)}{\tau^{\beta}}(\nabla U_h^{n+1-i},\nabla v_h)\\&+(\mathcal{F}(u_{H}^{n+1})+\mathcal{F}'(u_{H}^{n+1})(U_{h}^{n+1}-u_{H}^{n+1}),v_h)=(g^{n+1},v_h),
\end{split}\end{equation}
where $h\ll H$.
\begin{rem}
In the solving system above, we can seek a solution $u^{n+1}_H\in V_H$ on the coarse grid $\mathfrak{T}_H$ in the nonlinear system (\ref{3.32})-(\ref{3.33}), then get the solution $U^{n+1}_h\in V_h$ on the fine grid $\mathfrak{T}_h$ in the linear system (\ref{3.34})-(\ref{3.35}). We call the system (\ref{3.32})-(\ref{3.33}) with (\ref{3.34})-(\ref{3.35}) as two-grid FE system, which is more efficient than the standard nonlinear FE system (\ref{3.31})-(\ref{3.3}). In the results of numerical calculations, we will see the CUP-time used by two-grid FE scheme is less than that by standard nonlinear FE scheme.
\end{rem}
\par
In what follows, for the convenience of discussions on stability and a priori error analysis, we first give the following lemma.
\begin{lem}\label{lemma2}For series $\{w^{n}\}$, the following inequality holds
\begin{equation}\begin{split}\label{3.7} \Big{(}\frac32\delta^{n+1}_tw^{n+1}-\frac12\delta^{n}_tw^{n},w^{n+1}\Big{)}\geq&
\frac{1}{4\tau}[\Lambda(w^{n+1},w^{n})-\Lambda(w^{n},w^{n-1})],
\end{split}
\end{equation}
where
\begin{equation}\begin{split}\label{5.8}
&
\Lambda(w^{n},w^{n-1})\triangleq\|w^{n}\|^2+\|2w^{n}-w^{n-1}\|^2.
\end{split}
\end{equation}
\end{lem}
\par
In the next process,
we firstly consider the stability for systems (\ref{3.32})-(\ref{3.33}) and (\ref{3.34})-(\ref{3.35}).
\section{Analysis of stability based on two-grid FE algorithm }
\label{sec:3}
We first derive the stability based on two-grid FE algorithm.
\begin{thm}\label{thm11}For the two-grid FE system (\ref{3.32})-(\ref{3.35}) based on coarse grid $\mathfrak{T}_H$ and fine grid $\mathfrak{T}_h$, the following stable inequality for $U_h^n\in V_h$ holds
\begin{equation}\begin{split}\label{3.5}
&\|U_h^{n}\|^2\leq C(\|U_{h}^{0}\|^2+\|u_H^{0}\|^2
+\max_{0\leq i\leq n}\|g^{i}\|^2),
\end{split}
\end{equation}
\end{thm}
\textbf{Proof.} We first consider the results for the case $n\geq1$.
Setting $v_h=U_h^{n+1}$ in (\ref{3.35}),
and noting that the inequality (\ref{3.7}),
we have
\begin{equation}\begin{split}\label{3.8}
&\frac{1}{4\tau}[\Lambda(U_h^{n+1},U_h^{n})-\Lambda(U_h^{n},U_h^{n-1})]
\\&+\sum_{i=0}^{n+1}\frac{p_{\alpha}(i)}{\tau^{\alpha}}(U_h^{n+1-i},U_h^{n+1})
+\sum_{i=0}^{n+1}\frac{p_{\beta}(i)}{\tau^{\beta}}(\nabla U_h^{n+1-i},\nabla U_h^{n+1})
\\=&-(\mathcal{F}(u_{H}^{n+1})+\mathcal{F}'(u_{H}^{n+1})(U_{h}^{n+1}-u_{H}^{n+1}),U_h^{n+1})+(g^{n+1},U_h^{n+1}).
\end{split}\end{equation}
Using Cauchy-Schwarz inequality and Young inequality, we easily get
\begin{equation}\begin{split}\label{3.81}
&\frac{1}{4\tau}[\Lambda(U_h^{n+1},U_h^{n})-\Lambda(U_h^{n},U_h^{n-1})]
\\&+\sum_{i=0}^{n+1}\frac{p_{\alpha}(i)}{\tau^{\alpha}}(U_h^{n+1-i},U_h^{n+1})
+\sum_{i=0}^{n+1}\frac{p_{\beta}(i)}{\tau^{\beta}}(\nabla U_h^{n+1-i},\nabla U_h^{n+1})
\\\leq&C\|u_{H}^{n+1}\|\|U_h^{n+1}\|+\|\mathcal{F}'(u_{H}^{n+1})\|_{\infty}(\|U_{h}^{n+1}\|^2\\&+\|u_{H}^{n+1}\|\|U_h^{n+1}\|)+\|g^{n+1}\|\|U_h^{n+1}\|
\\\leq&C(\|u_{H}^{n+1}\|^2+\|U_h^{n+1}\|^2+\|g^{n+1}\|^2).
\end{split}\end{equation}
Sum (\ref{3.81}) for $n$ from $1$ to $L$ and use (\ref{3.7}) to get
\begin{equation}\begin{split}\label{3.82}
&\Lambda(U_h^{n+1},U_h^{n})+\tau^{1-\alpha}\sum_{n=1}^{L}\sum_{i=0}^{n+1}p_{\alpha}(i)(U_h^{n+1-i},U_h^{n+1})
\\&+\tau^{1-\beta}\sum_{n=1}^{L}\sum_{i=0}^{n+1}p_{\beta}(i)(\nabla U_h^{n+1-i},\nabla U_h^{n+1})\\\leq&C\tau\sum_{n=1}^{L}(\|u_{H}^{n+1}\|^2+\|U_h^{n+1}\|^2+\|g^{n+1}\|^2).
\end{split}\end{equation}
Set $v_h=U^1_h$ in (\ref{3.34}) and use Cauchy-Schwarz inequality and Young inequality to arrive at
\begin{equation}\begin{split}\label{3.83}
&\Big{(}\delta^{1}_t U^{1}_h,U^1_h\Big{)}+\sum_{i=0}^{1}\frac{p_{\alpha}(i)}{\tau^{\alpha}}(U_h^{1-i},U_h^1)+\sum_{i=0}^{1}\frac{p_{\beta}(i)}{\tau^{\beta}}(\nabla U_h^{1-i},\nabla U_h^1)\\=&-(\mathcal{F}(u_{H}^{1})+\mathcal{F}'((u_{H}^{1})(U_{h}^{1}-u_{H}^{1}),U_h^1)+(g^{1},U_h^1)
\\\leq&C(\|u_{H}^{1}\|^2+\|U_h^{1}\|^2+\|g^{1}\|^2).
\end{split}\end{equation}
From (\ref{3.83}), it easily follows that
\begin{equation}\begin{split}\label{3.85}
&\|U^{1}_h\|^2+\tau^{1-\alpha}\sum_{i=0}^{1}p_{\alpha}(i)(U_h^{1-i},U_h^1)+\tau^{1-\beta}\sum_{i=0}^{1}p_{\beta}(i)(\nabla U_h^{1-i},\nabla U_h^1)
\\\leq&C\tau(\|u_{H}^{1}\|^2+\|U_h^{0}\|^2+\|U_h^{1}\|^2+\|g^{1}\|^2).
\end{split}\end{equation}
Make a combination for (\ref{3.82}) and (\ref{3.85}) to get
\begin{equation}\begin{split}\label{3.84}
&\|U_h^{L}\|^2+\tau^{1-\alpha}\sum_{n=0}^{L}\sum_{i=0}^{n}\frac{p_{\alpha}(i)}{\tau^{\alpha}}(U_h^{n-i},U_h^{n})
+\tau^{1-\beta}\sum_{n=0}^{L}\sum_{i=0}^{n}\frac{p_{\beta}(i)}{\tau^{\beta}}(\nabla U_h^{n-i},\nabla U_h^{n})\\\leq&C\tau\sum_{n=0}^{L}(\|u_{H}^{n}\|^2+\|U_h^{n}\|^2+\|g^{n}\|^2)+C\|U_h^0\|^2.
\end{split}\end{equation}
Note that lemma \ref{lm1} and use Cronwall lemma to get
\begin{equation}\begin{split}\label{3.86}
&\|U_h^{L}\|^2\leq C\|U_h^0\|^2+C\tau\sum_{n=0}^{L}(\|u_{H}^{n}\|^2+\|g^{n}\|^2).
\end{split}\end{equation}
For the next estimates, we have to discuss the term $\|u_{H}^{n}\|^2$.
\par
In (\ref{3.32}) and (\ref{3.33}), we take $u^{1}_H$ and $u^{n+1}_H$ for $v_H$, respectively, and use a similar process of derivation to the $\|U_h^{L}\|^2$ to arrive at
\begin{equation}\begin{split}\label{3.87}
&\|u_H^{n}\|^2\leq C(\|u_H^{0}\|^2
+\max_{0\leq i\leq n}\|g^{i}\|^2)
\end{split}
\end{equation}
Substitute (\ref{3.87}) into (\ref{3.86}) and note that $\tau\sum_{n=0}^{L}\leq T$ to get
\begin{equation}\begin{split}\label{3.88}
&\|U_h^{L}\|^2\leq C(\|U_{h}^{0}\|^2+\|u_{H}^{0}\|^2+\max_{0\leq i\leq L}\|g^{i}\|^2),
\end{split}\end{equation}
which indicate that the conclusion (\ref{3.5}) of theorem \ref{thm11} holds.

\section{Error analysis based on two-grid algorithm}
\label{sec:4}
For discussing and deriving a priori error estimates based on fully discrete two-grid FE method, we have to introduce a Ritz-projection operator which is defined by finding $\Psi_\hslash: H_0^1(\Omega)\rightarrow V_\hslash$ such that
\be\label{3.05}(\nabla(\Psi_\hslash w),\nabla w_\hslash)=(\nabla w,\nabla w_\hslash),\forall w_\hslash\in V_\hslash,\ee
with the following estimate
\be\label{3.06}\|w-\Psi_\hslash w\|+\hslash\|w-\Psi_\hslash w\|_{1}\leq C\hslash^{r+1}\|w\|_{r+1},\forall w\in H_0^1(\Omega)\cap H^{r+1}(\Omega),\ee
where $\hslash$ is coarse grid step length $H$ or fine grid size $h$ and the norms are defined by $\|w\|_l=\Big{(}\sum\limits_{0\leq |\theta|\leq l}\int_{\Omega}|D^{\theta}w|^2d\textbf{x}\Big{)}^{\frac12}$ with the polynomial's degree $l$.
\par
In the following contents, based on the given Ritz-projection \cite{14} and estimate inequality (\ref{3.06}), we will do some detailed discussions on a priori error analysis.
\par
 Now we rewrite the errors as
$$u(t_n)-U^n_h=(u(t_n)-\Psi_h{U}^n)+(\Psi_h{U}^n-U^n_h)=\mathfrak{P}_u^n+\mathfrak{M}_u^n.$$
\begin{thm}\label{thm1.5}With $u(t_n)\in H_0^1(\Omega)\cap H^{r+1}(\Omega)$, $U_h^n\in V_h$ and $U_h^0=\Psi_h{u}(0)$, we obtain the following a priori error results in $L^2$-norm
\begin{equation}\begin{split}\label{5.11}
&\|u(t_n)-U^n_h\|\leq C[\tau^2+(1+\tau^{-\alpha})h^{r+1}+(1+\tau^{-2\alpha})H^{2r+2}],
\end{split}
\end{equation}
where $C$ is a positive constant independent of coarse grid step length $H$, fine grid size $h$ and time step parameters $\tau$.
\end{thm}
{\textbf{Proof.}}
Combine (\ref{3.35}) and (\ref{3.1}) with (\ref{3.05}) to arrive at the error equations for any $v_h\in V_h$ and $n\geq1$
\begin{equation}\begin{split}\label{5.31}
&\Big{(}\frac32\delta^{n+1}_t\mathfrak{M}_u^{n+1}-\frac12\delta^{n}_t\mathfrak{M}_u^{n},v_h\Big{)}
\\&+\sum_{i=0}^{n+1}\frac{p_{\alpha}(i)}{\tau^{\alpha}}(\mathfrak{M}_u^{n+1-i},v_h)
+\sum_{i=0}^{n+1}\frac{p_{\beta}(i)}{\tau^{\beta}}(\nabla \mathfrak{M}_u^{n+1-i},\nabla v_h)\\=&-\Big{(}\frac32\delta^{n+1}_t\mathfrak{P}_u^{n+1}-\frac12\delta^{n}_t\mathfrak{P}_u^{n},v_h\Big{)}
-\sum_{i=0}^{n+1}\frac{p_{\alpha}(i)}{\tau^{\alpha}}(\mathfrak{P}_u^{n+1-i},v_h)
\\&-(\mathcal{F}(u^{n+1})-\mathcal{F}(u_{H}^{n+1})+\mathcal{F}'(u_{H}^{n+1})(\mathfrak{M}_u^{n+1}+\mathfrak{P}_u^{n+1}-u^{n+1}+u_{H}^{n+1}),v_h)
\\&+(\bar{e}^{n+1}_1,v_h)+(\bar{e}^{n+1}_2,v_h)+(\Delta\bar{e}^{n+1}_3,v_h)
\\\doteq& I_1+I_2+I_3+I_4+I_5+I_6.
\end{split}
\end{equation}
In what follows, we need to estimate the terms $I_j,j=1,\cdots,6$. First we estimate the third term $I_3$.
For considering the nonlinear term, we use Taylor expansion to obtain
\begin{equation}\label{5.32}
\begin{aligned}
&\mathcal{F}(u^{n+1})-\mathcal{F}(u_{H}^{n+1})=\mathcal{F}'(u_{H}^{n+1})(u^{n+1}-u_{H}^{n+1})+\frac{1}{2}\mathcal{F}''(\chi^{n+1})(u^{n+1}-u_{H}^{n+1})^{2},
\end{aligned}
\end{equation}
where $\chi^{j}$ is a value between $u^{j}$ and $u_{H}^{j}$.
\\
Based on (\ref{5.32}), we obtain
\begin{equation}\begin{split}\label{3.13}
&\mathcal{F}(u^{n+1})-\mathcal{F}(u_{H}^{n+1})+\mathcal{F}'(u_{H}^{n+1})(\mathfrak{M}_u^{n+1}+\mathfrak{P}_u^{n+1}-u^{n+1}+u_{H}^{n+1})\\= &\mathcal{F}'(u_{H}^{n+1})(\mathfrak{M}_u^{n+1}+\mathfrak{P}_u^{n+1})+\frac{1}{2}\mathcal{F}''(\chi^{n+1})(u^{n+1}-u_{H}^{n+1})^{2}.
\end{split}
\end{equation}
So, we have
\begin{equation}\begin{split}\label{3.14}
&I_3=-(\mathcal{F}(u^{n+1})-\mathcal{F}(u_{H}^{n+1})+\mathcal{F}'(u_{H}^{n+1})(\mathfrak{M}_u^{n+1}+\mathfrak{P}_u^{n+1}-u^{n+1}+u_{H}^{n+1}),v_h)
\\\leq& \frac12\|\mathcal{F}'(u_{H}^{n+1})\|_{\infty}(\|\mathfrak{P}_u^{n+1}\|^2+\|\mathfrak{M}_u^{n+1}\|^2)
+\frac{1}{4}\|\mathcal{F}''(\chi^{n+1})\|_{\infty}\|(u^{n+1}-u_{H}^{n+1})^{2}\|^2
\\&+(\frac12\|\mathcal{F}'(u_{H}^{n+1})\|_{\infty}+\frac{1}{4}\|\mathcal{F}''(\chi^n)\|_{\infty})\|v_h\|^2.
\end{split}
\end{equation}
We now use Cauchy-Schwarz inequality with Young inequality to get
\begin{equation}\begin{split}\label{4.1}
I_1=&-\Big{(}\frac32\delta^{n+1}_t\mathfrak{P}_u^{n+1}-\frac12\delta^{n}_t\mathfrak{P}_u^{n},v_h\Big{)}
\\\leq &\Big{\|}\frac32\delta^{n+1}_t\mathfrak{P}_u^{n+1}-\frac12\delta^{n}_t\mathfrak{P}_u^{n}\Big{\|}\|v_h\|
\\\leq & C\int_{t_{n-1}}^{t_{n+1}}\|\mathfrak{P}_{ut}\|^2ds+C\|v_h\|^2,
\end{split}
\end{equation}
and
\begin{equation}\begin{split}\label{4.2}
I_4+I_5+I_6=&(\bar{e}^{n+1}_1,v_h)+(\bar{e}^{n+1}_2,v_h)+(\Delta\bar{e}^{n+1}_3,v_h)
\\\leq& C(\tau^4+\|v_h\|^2).
\end{split}
\end{equation}
By using lemma \ref{lemma4} with Cauchy-Schwarz inequality and Young inequality, we have
\begin{equation}\begin{split}\label{4.21}
I_2=&-\sum_{i=0}^{n+1}\frac{p_{\alpha}(i)}{\tau^{\alpha}}(\mathfrak{P}_u^{n+1-i},v_h)
\\\leq& \frac{1}{\tau^{\alpha}}\sum_{i=0}^{n+1}|p_{\alpha}(i)||(\mathfrak{P}_u^{n+1-i},v_h)|
\\\leq& \frac{\alpha+2}{2\tau^{\alpha}}g_{0}^{\alpha}\|\mathfrak{P}_u^{n+1}\|\|v_h\|
+\frac{1}{\tau^{\alpha}}\sum_{i=1}^{n+1}\Big{|}\frac{\alpha+2}{2}g_{i}^{\alpha}+\frac{-\alpha}{2}g_{i-1}^{\alpha}\Big{|}\|\mathfrak{P}_u^{n+1-i}\|\|v_h\|
\\\leq& Ch^{r+1}\|v_h\|\Big{(}\frac{\alpha+2}{2\tau^{\alpha}}g_{0}^{\alpha}
+\frac{1}{\tau^{\alpha}}\sum_{i=1}^{n+1}\Big{|}-\frac{\alpha+2}{2}g_{i}^{\alpha}\Big{|}
+\frac{1}{\tau^{\alpha}}\sum_{i=1}^{n+1}\Big{|}\frac{-\alpha}{2}g_{i-1}^{\alpha}\Big{|}\Big{)}
\\\leq& Ch^{r+1}\|v_h\|\Big{(}\frac{\alpha+1}{\tau^{\alpha}}g_{0}^{\alpha}
+\frac{\alpha+2}{2\tau^{\alpha}}\sum_{i=1}^{n+1}-g_{i}^{\alpha}
+\frac{\alpha}{2\tau^{\alpha}}\sum_{i=1}^{n}-g_{i}^{\alpha}\Big{)}
\\\leq& C(\alpha)\tau^{-\alpha}h^{r+1}\|v_h\|
\\\leq& C(\alpha)\tau^{-2\alpha}h^{2r+2}+C\|v_h\|^2.
\end{split}
\end{equation}
In (\ref{5.31}), (\ref{3.14})-(\ref{4.21}), we take $v_h=\mathfrak{M}_u^{n+1}$ and make a combination for these expressions to get
\begin{equation}\begin{split}\label{7.1}
&\frac{1}{4\tau}[\Lambda(\mathfrak{M}_u^{n+1},\mathfrak{M}_u^{n})-\Lambda(\mathfrak{M}_u^{n},\mathfrak{M}_u^{n-1})]
\\&+\sum_{i=0}^{n+1}\frac{p_{\alpha}(i)}{\tau^{\alpha}}(\mathfrak{M}_u^{n+1-i},\mathfrak{M}_u^{n+1})
+\sum_{i=0}^{n+1}\frac{p_{\beta}(i)}{\tau^{\beta}}(\nabla \mathfrak{M}_u^{n+1-i},\nabla \mathfrak{M}_u^{n+1})
\\\doteq& I_1+I_2+I_3+I_4+I_5+I_6
\\\leq& C(\tau^4+\tau^{-2\alpha}h^{2r+2})+\frac12\|\mathcal{F}'(u_{H}^{n+1})\|_{\infty}\|\mathfrak{P}_u^{n+1}\|^2
\\&+\frac{1}{4}\|\mathcal{F}''(\chi^{n+1})\|_{\infty}\|(u^{n+1}-u_{H}^{n+1})^{2}\|^2
\\&+(\|\mathcal{F}'(u_{H}^{n+1})\|_{\infty}+\frac{1}{4}\|\mathcal{F}''(\chi^n)\|_{\infty}+1)\|\mathfrak{M}_u^{n+1}\|^2.
\end{split}
\end{equation}
Multiply (\ref{7.1}) by $4\tau$ and sum (\ref{7.1}) for $n$ from $1$ to $L$ to get
\begin{equation}\begin{split}\label{4.25}
&\Lambda(\mathfrak{M}_u^{L+1},\mathfrak{M}_u^{L})
+4\tau^{1-\alpha}\sum_{n=1}^{L}\sum_{i=0}^{n+1}p_{\alpha}(i)(\mathfrak{M}_u^{n+1-i},\mathfrak{M}_u^{n+1})
\\&+4\tau^{1-\beta}\sum_{n=1}^{L}\sum_{i=0}^{n+1}p_{\beta}(i)(\nabla \mathfrak{M}_u^{n+1-i},\nabla \mathfrak{M}_u^{n+1})
\\\leq& \Lambda(\mathfrak{M}_u^{1},\mathfrak{M}_u^{0})+C\tau\sum_{n=1}^{L}(\tau^4+\tau^{-2\alpha}h^{2r+2})
+C\tau\sum_{n=1}^{L}\|\mathcal{F}'(u_{H}^{n+1})\|_{\infty}\|\mathfrak{P}_u^{n+1}\|^2
\\&+\tau\sum_{n=1}^{L}\|\mathcal{F}''(\chi^{n+1})\|_{\infty}\|(u^{n+1}-u_{H}^{n+1})^{2}\|^2
\\&+\tau\sum_{n=1}^{L}(\|\mathcal{F}'(u_{H}^{n+1})\|_{\infty}+\frac{1}{4}\|\mathcal{F}''(\chi^n)\|_{\infty}+1)\|\mathfrak{M}_u^{n+1}\|^2.
\end{split}
\end{equation}
Subtract (\ref{3.34}) from (\ref{3.30}), we have
\begin{equation}\begin{split}\label{3.39}
&\Big{(}\delta^{1}_t\mathfrak{M}_u^{1},v_h\Big{)}+\sum_{i=0}^{1}\frac{p_{\alpha}(i)}{\tau^{\alpha}}(\mathfrak{M}_u^{1-i},v_h)
+\sum_{i=0}^{1}\frac{p_{\beta}(i)}{\tau^{\beta}}(\nabla \mathfrak{M}_u^{1-i},\nabla v_h)\\=&-\Big{(}\delta^{n+1}_t\mathfrak{P}_u^{1},v_h\Big{)}
-\sum_{i=0}^{1}\frac{p_{\alpha}(i)}{\tau^{\alpha}}(\mathfrak{P}_u^{1-i},v_h)-(\mathcal{F}(u^{1})-\mathcal{F}(u_{H}^{1})\\&
-\mathcal{F}'(u_{H}^{1})(U_{h}^{1}-u_{H}^{1}),v_h)+(\bar{e}^1_1,v_h)+(\bar{e}^1_2,v_h)+(\Delta\bar{e}^1_3,v_h).
\end{split}\end{equation}
In (\ref{3.39}), we choose $v_h=\mathfrak{M}_u^{1}$ and use (\ref{3.14}) and (\ref{4.21}) to get
\begin{equation}\begin{split}\label{3.42}
&\|\mathfrak{M}_u^{1}\|^2-\|\mathfrak{M}_u^{0}\|^2+\|\mathfrak{M}_u^{1}-\mathfrak{M}_u^{0}\|^2
\\&+2\tau^{1-\alpha}\sum_{i=0}^{1}p_{\alpha}(i)(\mathfrak{M}_u^{1-i},\mathfrak{M}_u^{1})
+2\tau^{1-\beta}\sum_{i=0}^{1}p_{\beta}(i)(\nabla \mathfrak{M}_u^{1-i},\nabla \mathfrak{M}_u^{1})\\=&-2\tau\Big{(}\delta^{n+1}_t\mathfrak{P}_u^{1},\mathfrak{M}_u^{1}\Big{)}
-2\tau\sum_{i=0}^{1}\frac{p_{\alpha}(i)}{\tau^{\alpha}}(\mathfrak{P}_u^{1-i},\mathfrak{M}_u^{1})-2\tau(\mathcal{F}(u^{1})-\mathcal{F}(u_{H}^{1})\\&
-\mathcal{F}'(u_{H}^{1})(U_{h}^{1}-u_{H}^{1}),\mathfrak{M}_u^{1})
+2\tau(\bar{e}^1_1,\mathfrak{M}_u^{1})+2\tau(\bar{e}^1_2,\mathfrak{M}_u^{1})+2\tau(\Delta\bar{e}^1_3,\mathfrak{M}_u^{1})
\\\leq& \tau\|\mathcal{F}'(u_{H}^{1})\|_{\infty}(\|\mathfrak{P}_u^{1}\|^2+\|\mathfrak{M}_u^{1}\|^2)
+\frac{\tau}{2}\|\mathcal{F}''(\chi^{1})\|_{\infty}\|(u^{1}-u_{H}^{1})^{2}\|^2
\\&+(\tau\|\mathcal{F}'(u_{H}^{1})\|_{\infty}+\frac{\tau}{2}\|\mathcal{F}''(\chi^1)\|_{\infty})\|\mathfrak{M}_u^{1}\|^2+C\tau^4
+C\tau^{-2\alpha}h^{2r+2}+\frac14\|\mathfrak{M}_u^{1}\|^2,
\end{split}\end{equation}
Simplifying for (\ref{3.42}) and using triangle inequality, we have
\begin{equation}\begin{split}\label{3.40}
&\Lambda(\mathfrak{M}_u^{1},\mathfrak{M}_u^{0})+2\tau^{1-\alpha}\sum_{i=0}^{1}p_{\alpha}(i)(\mathfrak{M}_u^{1-i},\mathfrak{M}_u^{1})
+2\tau^{1-\beta}\sum_{i=0}^{1}p_{\beta}(i)(\nabla \mathfrak{M}_u^{1-i},\nabla \mathfrak{M}_u^{1})
\\\leq& C\tau h^{2r+2}+C\tau^{-2\alpha}h^{2r+2}
+\frac{\tau}{2}\|(u^{1}-u_{H}^{1})^{2}\|^2+C\tau^4.
\end{split}\end{equation}
Combine (\ref{4.25}) with (\ref{3.40}) and note that $\mathfrak{M}_u^{0}=0$ to get
\begin{equation}\begin{split}\label{4.26}
&\Lambda(\mathfrak{M}_u^{L+1},\mathfrak{M}_u^{L})
+4\tau^{1-\alpha}\sum_{n=-1}^{L}\sum_{i=0}^{n+1}p_{\alpha}(i)(\mathfrak{M}_u^{n+1-i},\mathfrak{M}_u^{n+1})
\\&+4\tau^{1-\beta}\sum_{n=-1}^{L}\sum_{i=0}^{n+1}p_{\beta}(i)(\nabla \mathfrak{M}_u^{n+1-i},\nabla \mathfrak{M}_u^{n+1})
\\\leq& C\tau\sum_{n=1}^{L}(\tau^4+\tau^{-2\alpha}h^{2r+2})
+C\tau\sum_{n=1}^{L}\|\mathcal{F}'(u_{H}^{n+1})\|_{\infty}\|\mathfrak{P}_u^{n+1}\|^2
\\&+\tau\sum_{n=1}^{L}\|\mathcal{F}''(\chi^{n+1})\|_{\infty}\|(u^{n+1}-u_{H}^{n+1})^{2}\|^2
+\tau\sum_{n=1}^{L}(\|\mathcal{F}'(u_{H}^{n+1})\|_{\infty}\\&+\frac{1}{4}\|\mathcal{F}''(\chi^n)\|_{\infty}+1)\|\mathfrak{M}_u^{n+1}\|^2+C\tau h^{2r+2}+
\frac{\tau}{2}\|(u^{1}-u_{H}^{1})^{2}\|^2+C\tau^4.
\end{split}
\end{equation}
By using the Cronwall lemma and the relationship (\ref{3.3}), we have for sufficiently small $\tau$
\begin{equation}\begin{split}\label{4.26}
\Lambda(\mathfrak{M}_u^{L+1},\mathfrak{M}_u^{L})\leq& C(\tau^4+\tau^{-2\alpha}h^{2r+2}+h^{2r+2})+C\tau\sum_{n=0}^{L}\|(u^{n+1}-u_{H}^{n+1})^{2}\|^2.
\end{split}
\end{equation}
For the next discussion, we need to give the estimate for the term $\|(u^{n+1}-u_{H}^{n+1})^{2}\|$.
\par
Subtract (\ref{3.32}), (\ref{3.33}) from (\ref{3.30}), (\ref{3.10}), respectively and use the Ritz-projection (\ref{3.05})
to arrive at the error equations under the coarse grid for any $v_H \in V_H$
\\
Case $n=0$:
\begin{equation}\begin{split}\label{5.30}
&\Big{(}\delta^{1}_t\mathfrak{D}_u^{1},v_H\Big{)}+\sum_{i=0}^{1}\frac{p_{\alpha}(i)}{\tau^{\alpha}}(\mathfrak{D}_u^{1-i},v_H)
+\sum_{i=0}^{1}\frac{p_{\beta}(i)}{\tau^{\beta}}(\nabla \mathfrak{D}_u^{1-i},\nabla v_H)\\=&
-\Big{(}\delta^{1}_t\mathfrak{A}_u^{1},v_H\Big{)}
-\sum_{i=0}^{1}\frac{p_{\alpha}(i)}{\tau^{\alpha}}(\mathfrak{A}_u^{1-i},v_H)
\\&-(\mathcal{F}(u^{1})-\mathcal{F}(u_H^{1}),v_H)+(\bar{e}^1_1,v_H)+(\bar{e}^1_2,v_H)+(\Delta\bar{e}^1_3,v_H),
\end{split}
\end{equation}
Case $n\geq 1$:
\begin{equation}\begin{split}\label{5.3}
&\Big{(}\frac32\delta^{n+1}_t\mathfrak{D}_u^{n+1}-\frac12\delta^{n}_t\mathfrak{D}_u^{n},v_H\Big{)}
\\&+\sum_{i=0}^{n+1}\frac{p_{\alpha}(i)}{\tau^{\alpha}}(\mathfrak{D}_u^{n+1-i},v_H)
+\sum_{i=0}^{n+1}\frac{p_{\beta}(i)}{\tau^{\beta}}(\nabla \mathfrak{D}_u^{n+1-i},\nabla v_H)\\=&-\Big{(}\frac32\delta^{n+1}_t\mathfrak{A}_u^{n+1}-\frac12\delta^{n}_t\mathfrak{A}_u^{n},v_H\Big{)}
-\sum_{i=0}^{n+1}\frac{p_{\alpha}(i)}{\tau^{\alpha}}(\mathfrak{A}_u^{n+1-i},v_H)
\\&-(\mathcal{F}(u^{n+1})-\mathcal{F}(u_{H}^{n+1}),v_H)+(\bar{e}^{n+1}_1,v_H)+(\bar{e}^{n+1}_2,v_H)+(\Delta\bar{e}^{n+1}_3,v_H),
\end{split}
\end{equation}
where $\mathfrak{A}_u^n=u(t_n)-\Psi_H{u}^n$, $\mathfrak{D}_u^n=\Psi_H{u}^n-u^n_H$.
\par
In (\ref{5.30}) and (\ref{5.3}), we take $v_H=\mathfrak{D}_u^{1}$ and $v_H=\mathfrak{D}_u^{n+1}$, respectively, and use a similar process of proof to the estimate for $\|u^n-U_h^n\|$ to get
\begin{equation}\begin{split}\label{5.33}
\|u^{n+1}-u_{H}^{n+1}\|\leq C(\tau^2+\tau^{-\alpha}H^{r+1}+H^{r+1}).
\end{split}
\end{equation}
Substitute the above estimate inequality (\ref{5.33}) into (\ref{4.26})
\begin{equation}\begin{split}\label{4.27}
\Lambda(\mathfrak{M}_u^{L+1},\mathfrak{M}_u^{L})\leq& C(\tau^4+\tau^{-2\alpha}h^{2r+2}+h^{2r+2}+\tau^{-4\alpha}H^{4r+4}+H^{4r+4}),
\end{split}
\end{equation}
which combine the triangle inequality with (\ref{3.06}) to get the conclusion of theorem \ref{thm1.5}.
\begin{rem}
Based on the theorem's results,
we obtain the temporal convergence rate with second-order result, which is free of fractional parameters $\alpha$ and $\beta$. Moreover, we can find that the time convergence rate by using second-order backward difference method and second-order WSGD scheme is higher than the one with $O(\tau^{2-\alpha}+\tau^{2-\beta})$ obtained by L1-approximation.
\end{rem}

\section{Numerical Tests}
\label{sec:5}
In this section, we consider a numerical example in space-time domain $[0,1]\times[0,1]^2$ to verify the theoretical results of two-grid FE algorithm combined with second-order backward difference method and second-order WSGD scheme. We now choose the nonlinear term $\mathcal{F}(u)=u^3-u$ and the exact solution $u(x,y,t)=t^2\sin(2\pi x)\sin(2\pi y)$, then easily determine that the known source function in (\ref{0.1}) is \begin{equation*}\begin{split}g(x,y,t)=\Big{[}2t-t^2+\frac{2t^{2-\alpha}}{\Gamma(3-\alpha)}&+16\pi^2\frac{t^{2-\beta}}{\Gamma(3-\beta)}\Big{]}\sin(2\pi x)\sin(2\pi y)\\&+t^6\sin^3(2\pi x)\sin^3(2\pi y).\end{split}\end{equation*}
\par
We now divide uniformly the spatial domain $[0,1]^2$ by using rectangular meshes, approximate first-order integer derivative with two-step backward Euler method and discretize the fractional direvative with second-order scheme. Now we take the continuous bilinear functions space $V_h$ with  $Q(x,y)=a_0+a_1x+a_2y+a_3xy$.
\par
For showing the current method in the this paper, we calculate some error results with convergence order for different fractional parameters $\alpha$ and $\beta$. In Table \ref{tab:1}, by taking fractional parameters $\alpha=0.01$, $\beta=0.99$ and fixed temperal step length $\tau=1/100$, we show some a priori errors in $L^2$-norm and convergence orders for two-grid algorithm with coarse and fine meshes $H=\sqrt{h}=1/4,1/5,1/6,1/7$ and FE method with $h=1/16,1/25,1/36,1/49$. From Table \ref{tab:1}, ones can see that the results with second-order convergence rate by using our method is stable and the CPU-time in seconds for two-grid FE method is less than that by making use of the standard FE method. In Table \ref{tab:2}, we use the same computing method and spatial meshes as in Table \ref{tab:1}, then obtain the errors and convergence rates when taking $\alpha=0.5$, $\beta=0.5$ and $\tau=1/100$. The similar calculated results with $\alpha=0.99$, $\beta=0.01$ and $\tau=1/100$ are also shown in Table \ref{tab:3}. These numerical results shown in Tables \ref{tab:1}-\ref{tab:3} also tell ones that compared to widely used L1-formula with $O(\tau^{2-\alpha}+\tau^{2-\beta})$, the WSGD scheme in this paper can get second-order convergence rate. For obtaining the second-order accuracy in time, we here use the second-order backward difference method to approximate time direction. Compared to commonly used one-step backward Euler difference method, for getting the same calculated accuracy, the second-order backward difference method can reduce the number of iterations in time and save the calculating time. From the numerical results presented in Tables \ref{tab:1}-\ref{tab:2}, ones can see that the convergence order of our method is slightly higher than that of the standard nonlinear FE method, while in Table 3, our method shows the similar convergence rate to that of nonlinear FE method. These phenomena indicate that compared to nonlinear FE method, our method has advantages in solving the time fractional Cable equation covering the space-time partial derivative term $^R_0\partial_{t}^\beta\Delta u$ with larger fractional parameter $\beta\in (0,1)$.
\par
\begin{figure*}
  \includegraphics[width=0.75\textwidth]{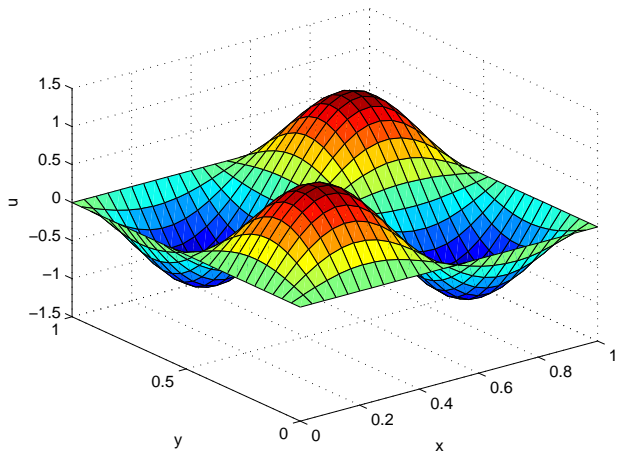}
\caption{Exact solution $u$}
\label{fig:1}       
\end{figure*}
\begin{figure*}
  \includegraphics[width=0.75\textwidth]{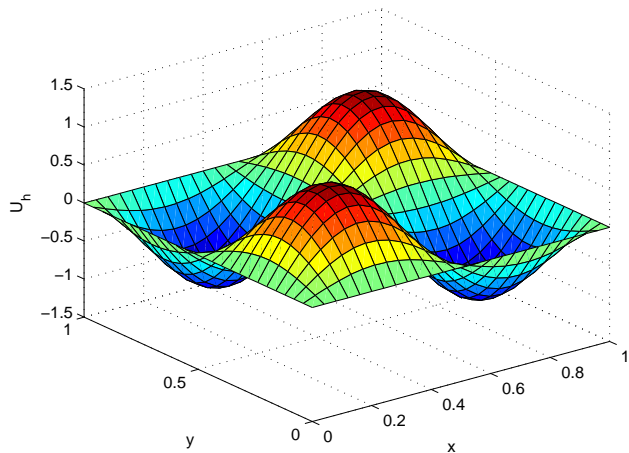}
\caption{Two-grid FE solution $U_h$}
\label{fig:2}
\end{figure*}
\begin{figure*}
  \includegraphics[width=0.75\textwidth]{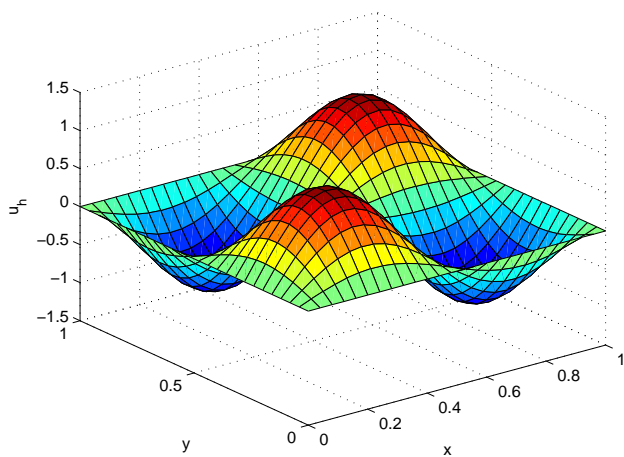}
\caption{FE solution $u_h$}
\label{fig:3}
\end{figure*}
\begin{figure*}
  \includegraphics[width=0.75\textwidth]{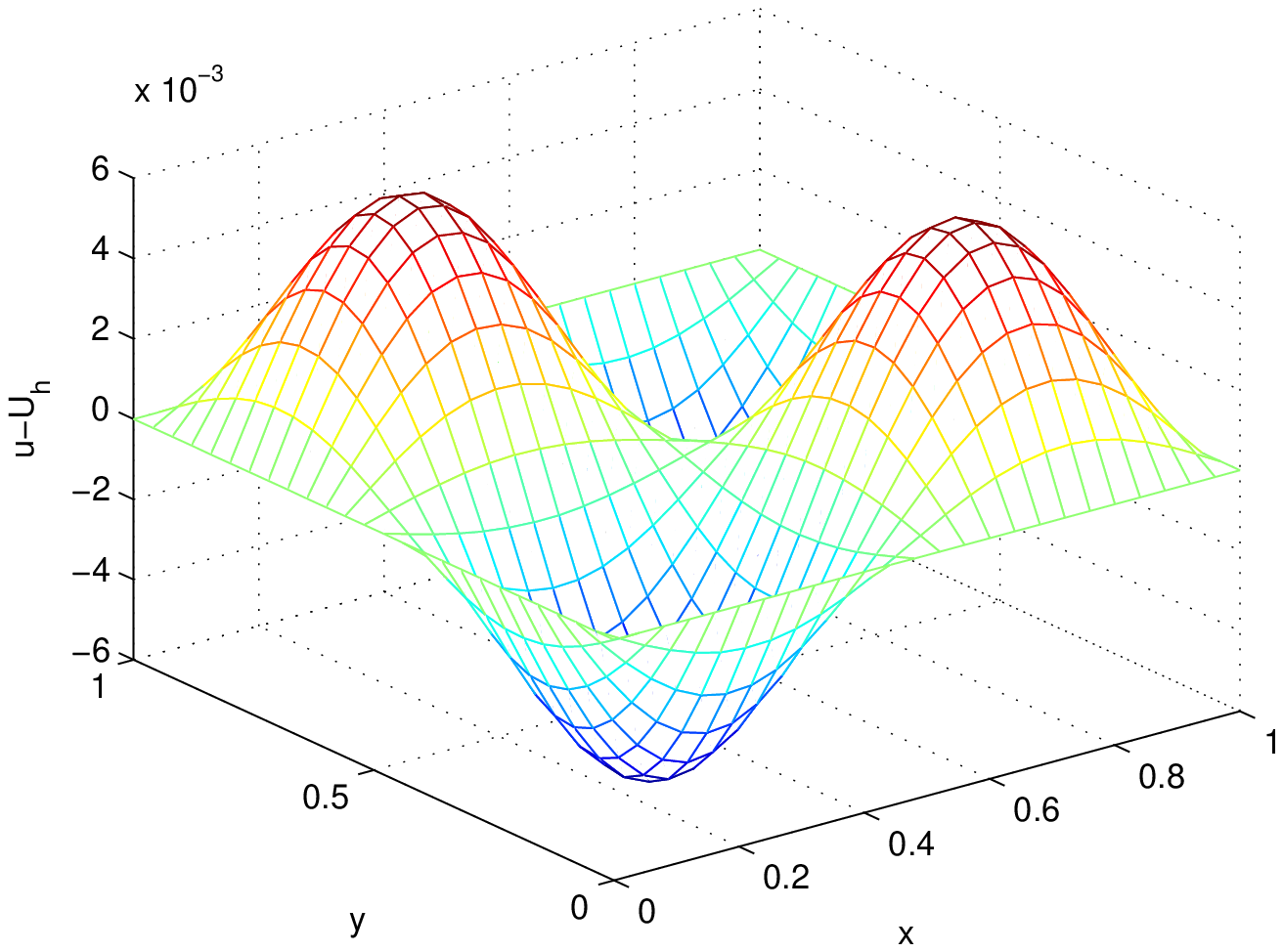}
\caption{The value for $u-U_h$}
\label{fig:4}
\end{figure*}
\begin{figure*}
  \includegraphics[width=0.75\textwidth]{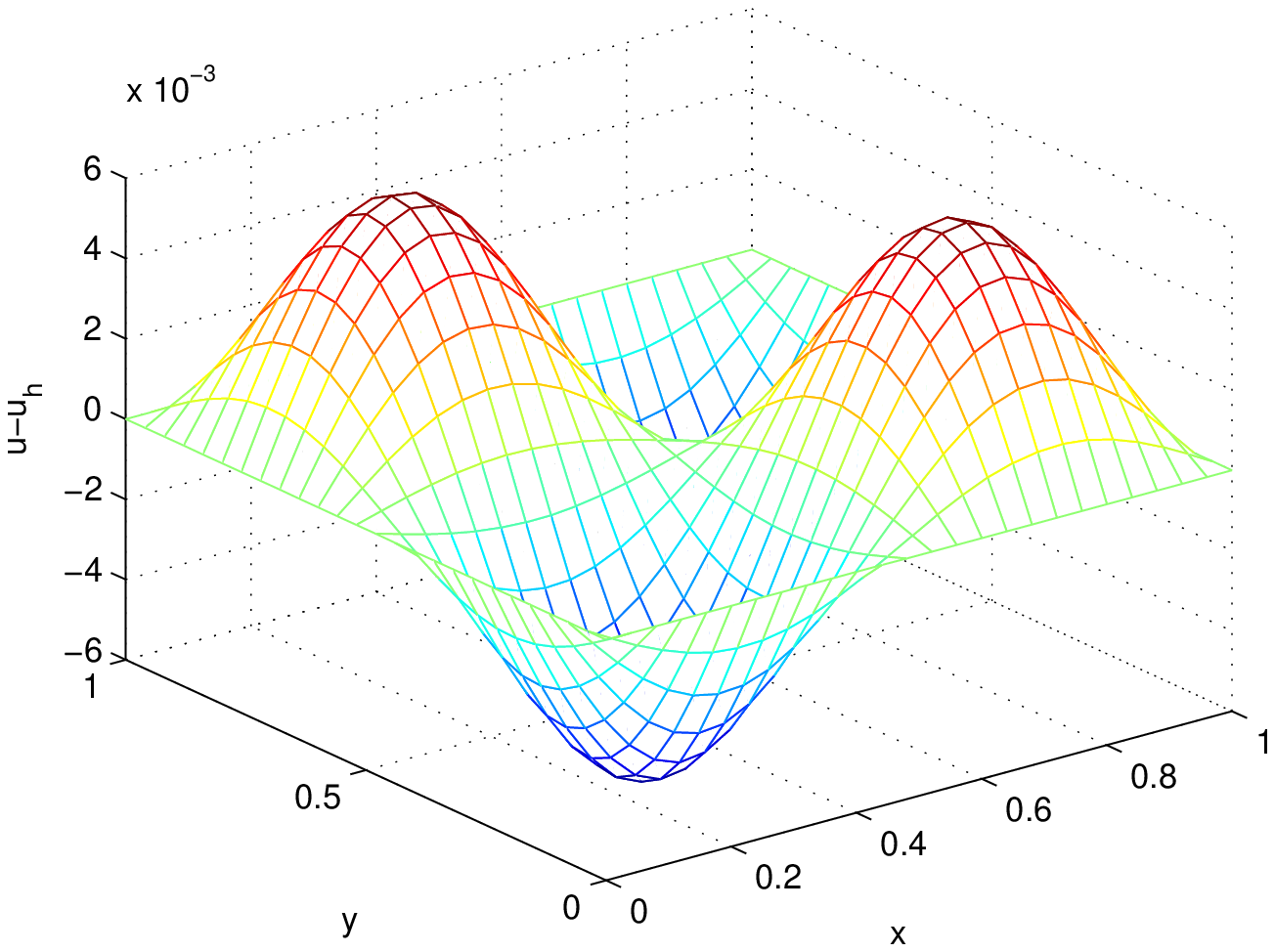}
\caption{The value for $u-u_h$}
\label{fig:5}
\end{figure*}

In Figures \ref{fig:1}-\ref{fig:3}, by taking $\alpha=0.99$, $\beta=0.01$, $\tau=1/100$ and $h=H^2=1/25$, we show the surfaces for the exact solution $u$, two-grid FE solution $U_h$ and FE solution $u_h$, respectively. Ones easily see that both two-grid FE solution $U_h$ based on coarse and fine meshes and FE solution $u_h$ can approximate well the exact solution $u$. Especially, from the surface for errors $u-U_h$ and $u-u_h$ in Figures \ref{fig:4}-\ref{fig:5}, ones easily find that two-grid FE method holds the similar computational accuracy to that of standard nonlinear FE method.
\par
In summary, from the computed error results and convergence rate in Tables \ref{tab:1}-\ref{tab:3} and the surfaces shown in Figures \ref{fig:1}-\ref{fig:5}, ones can know that with the similar computational accuracy to that for standard nonlinear FE method, our two-grid FE method is more efficient in computational time than the standard nonlinear FE method. Moreover, the current method combined with the second-order backward difference method and second-order WSGD operater can get a stable second-order convergence rate, which is independent of
fractional parameters $\alpha$ and $\beta$ and is higher than the convergence result $O(\tau^{2-\alpha}+\tau^{2-\beta})$ derived by L1-approximation.

\section{Some concluding remarks}
\label{sec:6}
In this article, we consider two-grid method combined with FE methods to give the numerical solution for nonlinear fractional Cable equations. First, we give some lemmas used in our paper; Second, we give the approximate formula for fractional derivative, then formulate the numerical scheme based on two-grid FE method; Finally, we do some detailed derivations for the stability of numerical scheme and a priori error analysis with second-order convergence rate in time, then compute some numerical errors and convergence orders to verify the theoretical results.
\par
From the numerical results, ones easily see that two-grid FE method studied in this paper can solve well the nonlinear time fractional Cable equation. Based on the point of view of calculating efficiency, compared to FE method, two-grid FE method can spend less time. Moreover, compared with the time convergence rate $O(\tau^{2-\alpha}+\tau^{2-\beta})$ obtained by usual L1-approximation, the current numerical scheme can arrive at second-order convergence rate independent of
fractional parameters $\alpha$ and $\beta$. Considering the mentioned advantages, in the future works, we will discuss the numerical theories of two-grid FE method for some space and space-time fractional partial differential equations with nonlinear term.

\section*{Acknowledgements} Authors thank the anonymous referees and editors very much for their valuable comments and suggestions,
which greatly help us improve this work. This work is supported by supported by the National Natural Science Fund (11301258, 11361035), Natural Science Fund of Inner Mongolia Autonomous Region (2015MS0114), and the Postgraduate Scientific Research Innovation Foundation of Inner Mongolia (No. 1402020201337).

%

%
\begin{table}
\caption{The $L^2$-errors with $\alpha=0.01$, $\beta=0.99$ and $\tau=1/100$}
\label{tab:1}
\begin{tabular}{llllll}
\hline\noalign{\smallskip}
  $H$ & $h$ & $\|u-U_{h}\|$ & Order & CPU time (in seconds)  \\
\noalign{\smallskip}\hline\noalign{\smallskip}
     1/4 & 1/16 & 6.3566e-003 & - & 43.231369  \\
     1/5 & 1/25 & 2.6118e-003 & 1.9930 &  122.612277 \\
     1/6 & 1/36 & 1.2323e-003 & 2.0600 &  414.219975 \\
     1/7 & 1/49 & 6.3532e-004 & 2.1489 &  1810.428349 \\
\noalign{\smallskip}\hline\noalign{\smallskip}
  FE algorithm & $h$ & $\|u-u_{h}\|$ & Order & CPU time (in seconds)  \\
\noalign{\smallskip}\hline\noalign{\smallskip}
      & 1/16 & 6.4246e-003 & - & 49.760083  \\
      & 1/25 & 2.6815e-003 & 1.9578 &  145.938800 \\
      & 1/36 & 1.3025e-003 & 1.9803 &  488.617402 \\
      & 1/49 & 7.0575e-004 & 1.9876 &  2112.565940 \\
\noalign{\smallskip}\hline
\end{tabular}
\end{table}

\begin{table}
\caption{The $L^2$-errors with $\alpha=0.5$, $\beta=0.5$ and $\tau=1/100$}
\label{tab:2}
\begin{tabular}{llllll}
\hline\noalign{\smallskip}
  $H$ & $h$ & $\|u-U_{h}\|$ & Order & CPU time (in seconds)  \\
\noalign{\smallskip}\hline\noalign{\smallskip}
     1/4 & 1/16 & 6.6252e-003 & - & 34.637650  \\
     1/5 & 1/25 & 2.7694e-003 & 1.9545 &  102.802120 \\
     1/6 & 1/36 & 1.3488e-003 & 1.9729 &  374.852060 \\
     1/7 & 1/49 & 7.3406e-004 & 1.9733 &  1745.224030 \\
 \noalign{\smallskip}\hline\noalign{\smallskip}
  FE algorithm & $h$ & $\|u-u_{h}\|$ & Order & CPU time (in seconds)  \\
\noalign{\smallskip}\hline\noalign{\smallskip}
      & 1/16 & 6.6292e-003 & - & 37.660242  \\
      & 1/25 & 2.7735e-003 & 1.9525 &  116.171144 \\
      & 1/36 & 1.3529e-003 & 1.9687 &  448.280514 \\
      & 1/49 & 7.3816e-004 & 1.9651 & 2189.299786  \\
  \noalign{\smallskip}\hline
\end{tabular}
\end{table}

\begin{table}
\caption{ The $L^2$-errors with $\alpha=0.99$, $\beta=0.01$ and $\tau=1/100$}
\label{tab:3}
\begin{tabular}{llllll}
\hline\noalign{\smallskip}
  $H$ & $h$ & $\|u-U_{h}\|$ & Order & CPU time (in seconds)  \\
 \noalign{\smallskip}\hline\noalign{\smallskip}
     1/4 & 1/16 & 6.9107e-003 & - & 53.860030  \\
     1/5 & 1/25 & 2.8841e-003 & 1.9581 &  153.053853 \\
     1/6 & 1/36 & 1.4003e-003 & 1.9815 &  497.518498 \\
     1/7 & 1/49 & 7.5807e-004 & 1.9905 &  2040.087416 \\
 \noalign{\smallskip}\hline\noalign{\smallskip}
  FE algorithm & $h$ & $\|u-u_{h}\|$ & Order & CPU time(in seconds)  \\
  \noalign{\smallskip}\hline\noalign{\smallskip}
      & 1/16 & 6.9107e-003 & - & 57.108134  \\
      & 1/25 & 2.8841e-003 & 1.9581 &  163.841391 \\
      & 1/36 & 1.4003e-003 & 1.9815 &  554.982539 \\
      & 1/49 & 7.5809e-004 & 1.9904 &  2416.876270 \\
  \noalign{\smallskip}\hline
\end{tabular}
\end{table}



\end{document}